\documentclass[12pt]{article}
\usepackage{graphicx}
\usepackage{amssymb, cite}
\usepackage{setspace} 
\usepackage[hmargin=1.15in,vmargin=1.15in]{geometry} 

\newcommand{\be}{\begin{equation}}
\newcommand{\ee}{\end{equation}}
\newcommand{\bqn}{\begin{eqnarray}}

\newcommand{\eqn}{\end{eqnarray}}

\newcommand{\bd}{\begin{description}}
\newcommand{\ed}{\end{description}}

\newtheorem{stat}{}[section]

\def\bs{\begin{stat}}
\def\es{\end{stat}}
\def\ben{\begin{enumerate}}
\def\een{\end{enumerate}}

\def\bp{\noindent{\bf Proof}  \ \ \ }
\newcommand{\ep}{\hfill $\square$}

\makeatletter
\@addtoreset{equation}{section}

\makeatother

\begin{document}

\begin{center}
{\large {\bf PACKING  3-VERTEX PATHS  
IN CLAW-FREE GRAPHS}}
\\[4ex]
{\large {\bf Alexander Kelmans}}
\\[2ex]
{\bf University of Puerto Rico, San Juan, Puerto Rico}
\\[0.5ex]
{\bf Rutgers University, New Brunswick, New Jersey}
\\[2ex]
\end{center}

\begin{abstract}
A $\Lambda $-{\em factor} of a graph $G$ is a spanning 
subgraph of $G$ whose every component is a 3-vertex 
path. 
Let $v(G)$  be the number 
of vertices of $G$.
A graph is {\em claw-free} if it does not have a subgraph isomorphic to $K_{1,3}$.
Our results include the following.
Let $G$ be a 3-connected claw-free graph, 
$x \in V(G)$, $e = xy \in E(G)$, and $L$ a 3-vertex path in $G$. Then 
$(c1)$ if $v(G) = 0 \bmod 3$, then $G$ has  a 
$\Lambda $-factor  containing (avoiding) $e$,
$(c2)$ if $v(G) = 1 \bmod 3$,  then $G - x$ has 
a $\Lambda $-factor, 
$(c3)$ if $v(G) = 2 \bmod 3$, then $G - \{x,y\}$ has a $\Lambda $-factor,
$(c4)$ if $v(G) = 0 \bmod 3$ and $G$ is either cubic or 4-connected, then $G - L$ has  a $\Lambda $-factor, and
$(c5)$ if $G$ is cubic and $E$ is a set of three edges in $G$, 
then $ G - E$ has  a $\Lambda $-factor if and only if the subgraph induced by $E$ in $G$ is not a claw and not a triangle.  
\\[1ex]
  
\end{abstract}

\section{Introduction}

\indent

We consider undirected graphs with no loops and 
no parallel edges. All notions and facts on graphs, that are  
used but not described here, can be found in \cite{BM,D,Wst}.
\\[1ex]
\indent
Given graphs $G$ and $H$, 
an $H$-{\em packing} of $G$ is a subgraph of $G$ 
whose every component is isomorphic to $H$.
An $H$-{\em packing} $P$ of $G$ is called 
an $H$-{\em factor} if $V(P) = V(G)$. 
The $H$-{\em packing problem}, i.e. the problem of 
finding in $G$ an $H$-packing, having the maximum 
number of vertices, turns out to be $NP$-hard if $H$ is 
a connected graph with at least three vertices \cite{HK}.
Let $\Lambda $ denote a 3-vertex path.
In particular, the $\Lambda $-packing problem
is $NP$-hard. Moreover, this problem remains 
$NP$-hard even for cubic graphs \cite{K1}.

Although the $\Lambda $-packing problem is $NP$-hard, 
i.e. possibly intractable in general, this problem turns out 
to be tractable for some natural classes of graphs
(see, for example, {\bf \ref{2conclfr}} below). 
It would be  also interesting to find polynomial-time algorithms 
that would provide a good approximation solution for 
the problem (e.g. {\bf \ref{km}}
and {\bf \ref{eb(G)clfr}} below).
In each case the corresponding packing problem is 
polynomially solvable.

Let $v(G)$ and $\lambda (G)$ denote the number 
of vertices and the maximum number of disjoint
3--vertex paths in $G$, respectively.
Obviously $\lambda (G) \le \lfloor v(G)/3 \rfloor $.

In \cite{K,KM} we answered the following natural question:
\\[1ex]
\indent
{\em How many disjoint 3-vertex paths must a cubic $n$-vertex graph have?}
\bs 
\label{km} 
If $G$ is a cubic graph, then 
$\lambda (G) \ge \lceil v(G)/4  \rceil$.
Moreover, there is a polynomial time algorithm for 
finding a $\Lambda $-packing having at least  
$\lceil v(G)/4  \rceil$ components.
\es

Obviously if every component of $G$ is $K_4$, then
$\lambda (G) = v(G)/4$. 
Therefore the bound in {\bf \ref{km}} is sharp.
\\[.5ex]
\indent
Let ${\cal G}^3_2 $ denote the set of graphs with each vertex of degree  $2$ or $3$.
In \cite{K} we  answered (in particular) the following  question:
\\[1ex]
\indent
{\em How many disjoint 3--vertex paths must an $n$-vertex  graph from ${\cal G}^3_2$ have?}

\bs 
\label{2,3-graphs} Suppose that $G \in {\cal G}^3_2$ and  
$G$ has no 5-vertex components.
Then $\lambda (G) \ge v(G)/4$.
\es

Obviously {\bf \ref{km}}  follows from {\bf \ref{2,3-graphs}} because if $G$ is a cubic graph, then $G \in {\cal G}^3_2$ and $G$ has no 5-vertex components. 
\\[1ex]
\indent
In \cite{K} we  also gave a construction  that allowed  to prove the following:
\bs 
\label{extrgraphs1}
There are infinitely many connected graphs for  which the bound  in {\bf \ref{2,3-graphs}} is attained.
Moreover, there are infinitely many subdivisions of
cubic 3-connected graphs for which the bound  in 
{\bf \ref{2,3-graphs}} is attained.
\es

The next interesting question is:
\\[1ex]
\indent
{\em How many disjoint 3-vertex paths must a cubic connected graph have?}
\\[1ex]
\indent
In \cite{K2} we proved the following.
Let ${\cal C}_n$ denote the set of connected cubic graphs with $n$ vertices.
\bs 
\label{cubic-connected} 
Let 
$\lambda _n = \min \{\lambda (G)/v(G): G \in {\cal C}_n\}$.
Then for some $c > 0$,
%
%
\[\frac{3}{11}(1 - \frac{c}{n}) \le \lambda _n \le 
\frac{3}{11}(1 - \frac{1}{n^2}).\]
\es

The similar question for cubic 2-connected graphs is still open:

\bs {\bf Problem.}
\label{2con-must}
How many disjoint 3-vertex paths must a cubic 2-connected graph have?
\es

It is known that

\bs 
\label{2con-nofactor}
There are infinitely many 2-connected and cubic graphs $G$ such that $\lambda (G) <   \lfloor v(G)/3 \rfloor $.
\es

Some such graph sequences were constructed in \cite{Kcntrex} to provide 2-connected counterexamples to Reed's domination conjecture.
Reed's conjecture claims that if $G$ is a connected cubic graph, then $\gamma (G) \le \lceil v(G)/3 \rceil $, where $\gamma (G)$ is the dominating number of $G$ (i.e. the size of a minimum vertex subset $X$ in $G$ such that every vertex in $G - X$ is adjacent to a vertex in $X$).
In particular, a graph sequence $(R_k:  k \ge 3)$ in \cite{Kcntrex} is such that each $R_k$ is a cubic graph of connectivity two and
$\gamma (G)/v(G) = \frac{1}{3} + \frac{1}{60}$.
 Obviously, $\gamma (G) \le v(G) - 2 \lambda (G)$.
 Therefore $\lambda (R_k)/ v(R_k) \le  \frac{13}{40}$.

The questions arise whether the claim of {\bf \ref{2con-nofactor}} is true for cubic  
2-connected  graphs having some additional properties. 
For example,

\bs {\bf Problem.} 
\label{Pr2con,cub,bip,pl}  
Is $\lambda (G) =   \lfloor v(G)/3 \rfloor $ true for every 
2-connected, cubic, bipartite, and planar  graph ?
\es

In \cite{K2con-cbp} we answered the question in 
{\bf \ref{Pr2con,cub,bip,pl}} by 
giving a construction that provides infinitely many  
 2-connected, cubic, bipartite, and planar  graphs such that  
$\lambda (G) <   \lfloor v(G)/3 \rfloor $.
\\[1ex]
\indent
As to cubic 3-connected graphs, an old open question here is:

\bs {\bf Problem.}
\label{Pr3con} Is the following claim true ?
\\[.5ex]
If $G$ is a 3-connected and cubic graph,
then $\lambda (G) =   \lfloor v(G)/3 \rfloor $.
\es

In \cite{K3con-cub} we discuss Problem {\bf \ref{Pr3con}} and
show, in particular, that the claim in 
{\bf \ref{Pr3con}} is equivalent to some seemingly 
much stronger claims. Here are some results of this kind.

\bs {\em \cite{K3con-cub}}
\label{cubic3-con}
The following are equivalent for cubic 3-connected graphs $G$:
\\[1ex]
$(z1)$
$v(G) = 0 \bmod 6$ $\Rightarrow$ $G$ has 
a $\Lambda $-factor,
\\[1ex] 
 $(z2)$ 
$v(G) = 0 \bmod 6$ $\Rightarrow$ for every 
$e \in E(G)$ there is a $\Lambda $-factor of $G$ 
avoiding $e$,
 \\[1ex]
 $(z3)$ 
$v(G) = 0 \bmod 6$ $\Rightarrow$ for every 
$e \in E(G)$ there is a $\Lambda $-factor of $G$ 
containing $e$,
 \\[1ex]
${(z4)}$
 $v(G) = 0 \bmod 6$ $\Rightarrow$ $G - X$ has 
a $\Lambda $-factor for every $X \subseteq E(G)$, 
$|X| = 2$,
\\[1ex]  
${(z5)}$
 $v(G) = 0 \bmod 6$ $\Rightarrow$ 
$G - L$ has a $\Lambda $-factor for every
3-vertex path $L$ in $G$,
 \\[1ex]
 $(t1)$ 
$v(G) = 2 \bmod 6$ $\Rightarrow$ $G - \{x,y\}$ 
has a $\Lambda $-factor for every $xy \in E(G)$,
 \\[1ex]
 $(f1)$ 
$v(G) = 4 \bmod 6$ $\Rightarrow$ $G - x$ 
has a $\Lambda $-factor for every $x \in V(G)$,
\\[1ex]
 $(f2)$ 
$v(G) = 4 \bmod 6$ $\Rightarrow$ $G - \{x, e\}$ 
has a $\Lambda $-factor for every $x \in V(G)$ and 
$e \in E(G)$.
\es

From {\bf \ref{cubic3-con}} it follows that if the claim in 
Problem {\bf \ref{Pr3con}} is true, then Reed's domination conjecture is true for 3-connected cubic graphs.
\\[1ex]
\indent
There are some interesting results on the $\Lambda $-packing problem for so called claw-free graphs.
A graph is called {\em claw-free} if it contains no induced subgraph isomorphic
to $K_{1,3}$ (which is called a {\em claw}). 

A vertex $x$ of a block $B$ in $G$ is called a {\em boundary vertex of} $B$ if $x$  belongs to another block of $G$.
If $B$ has exactly one boundary vertex,  then $B$ is called an {\em end-block} of $G$.
Let $eb(G)$ denote the number of end-blocks of $G$.
\bs {\em \cite{KKN}}
\label{2conclfr} 
Suppose that $G$ is a  claw-free graph and $G$ is either 
2-connected or connected with exactly two end-blocks.
Then $\lambda(G) = \lfloor v(G)/3 \rfloor$.
\es 

\bs {\em \cite{KKN}}
\label{eb(G)clfr}
Suppose that $G$ is a connected claw-free graph and  
$eb(G) \ge 2$. 
Then  $\lambda(G) \ge \lfloor (v(G) - eb(G) + 2)/3 \rfloor$,
and this lower bound is sharp.
\es

Obviously the claim in {\bf \ref{2conclfr}} about connected claw-free graphs with exactly two end-blocks follows from 
{\bf \ref{eb(G)clfr}}.
\\[1ex]
\indent
In this paper (see Section \ref{clfree}) we give some more results on the $\Lambda $-packings in 
claw-free graphs.
We show, in particular, the following: 
\\[1ex]
$(c1)$ all claims in  {\bf \ref{cubic3-con}} except for $(z5)$ are true for 3-connected claw-free graphs and $(z5)$ is true for 
cubic, 2-connected, and claw-free graphs distinct from $K_4$ (see {\bf \ref{clfree,3con,cubic,ztf}} below),
%
\\[1ex]
$(c2)$ if $G$ is a 3-connected  claw-free graph and
$v(G) = 0 \bmod 3$, then for every edge $e$ in $G$ there exists a $\Lambda $-factor of $G$ containing $e$ 
(see {\bf \ref{clfree,3con,einL}}),
\\[1ex]
$(c3)$ 
if $G$ is a 2-connected  claw-free graph and
$v(G) = 0 \bmod 3$, then for every edge  $e$ in $G$
there exists a $\Lambda $-factor of $G$ avoiding $e$, i.e. 
$G - e$ has a $\Lambda$-factor
(see {\bf \ref{clfree,2con-avoid-e}}),
%
\\[1ex]
$(c4)$ 
if $G$ is a 2-connected  claw-free graph and
$v(G) = 1 \bmod 3$, then $G - x$ has 
a $\Lambda$-factor for every vertex $x$ in $G$ 
(see {\bf \ref{clfree,2con,x-}}),
\\[1ex]
$(c5)$ if $G$ is a 3-connected  claw-free graph and
$v(G) = 2 \bmod 3$, then $G - \{x,y\}$ has 
a $\Lambda$-factor for every edge $xy$ in $G$
(see {\bf \ref{clfree,3con,(x,y)-}}),
\\[1ex]
$(c6)$ if $G$ is a 3-connected, claw-free, and  cubic graph with $v(G) \ge  6$ or a 4-connected claw-free graph, then for every 3-vertex path  $L$ in $G$ there exists a $\Lambda $-factor containing $L$, i.e. $G - L$ has a  $\Lambda $-factor
(see  {\bf \ref{2con,cubic,tr}} and 
{\bf \ref{clfree,4con,deg3cntrL-}}),
\\[1ex]
$(c7)$ if $G$ is a cubic, 3-connected, and  claw-free graph with $v(G) \ge 6$ and 
$E$ is a set of three edges in $G$, 
then $ G - E$ has  a $\Lambda $-factor if and only if the subgraph induced by $E$ in $G$ is not a claw and not a triangle 
(see {\bf \ref{clfree,2con-avoid(a,b,c)}}),
\\[1ex]
$(c8)$
if $G$ is a 3-connected  claw-free 
graph,  $v(G) = 1 \bmod 3$, $x \in V(G)$, and $e \in E(G)$,
then $G - \{x,e\}$ has a $\Lambda $-factor
(see {\bf \ref{clfree,3con,(x,e)-}}).

\section
{Main results}
\label{clfree}

\indent

Theorems  {\bf \ref{2conclfr}} and  {\bf \ref{eb(G)clfr}} above
describe some properties of maximum 
$\Lambda $-packings in claw-free graphs.
In this section we establish some more properties of
$\Lambda $-packings in claw-free graphs.

Let $G$ be a  graph and $B$ be a block of $G$, and so $B$ is either 2-connected or consists of two vertices and one edge.
As above, a vertex $x$ of $B$ is called a {\em boundary vertex of} $B$
if $x$  belongs to another block of $G$, and an {\em inner vertex of} $B$, otherwise. 
%
If $B$ has exactly one boundary vertex,  then $B$ is called an {\em end-block} of $G$.

Let $F$ be a graph, $x \in V(F)$, and $X = \{x_1,x_2, x_3\}$ be the set of  vertices in $F$ adjacent to $x$.
Let $T$ be a triangle, $V(T) = \{t_1, t_2, t_3\}$, and $V(F) \cap V(T) = \emptyset $. 
Let $G = (F - x) \cup T \cup \{x_it_i: i \in \{1,2,3\}$.
We say that $G$ {\em is obtained from $F$ by replacing a vertex $x$ by a triangle}. Let $F^\Delta $ denote the graph obtained from a cubic graph $F$ by replacing each 
vertex  of $F$ by a triangle. Obviously, $F^\Delta $ is claw-free, every vertex  belongs to exactly one triangle, and every edge belongs to at most one triangle in $F^\Delta $.
\bs
\label{2con,cubic,tr}
Let $G'$ be a cubic 2-connected graph and
$G $ be the graph obtained from $G'$ by replacing each 
vertex $v$  of $G'$ by a triangle $\Delta _v$.
 Let $L$ be a 3-vertex path in $G$.
Then
\\[0.5ex]
$(a)$ $G - L$ has a $\Lambda$-factor.

Moreover,
\\[0.5ex]
$(a1)$ if $L$ induces  a triangle in $G$, then 
$G$ has a $\Lambda$-factor $R$ containing $L$ and such that each component of $R$ induces a triangle
\\[0.5ex]
$(a2)$ if $L$ does not induce a triangle in $G$, then 
$G$ has a $\Lambda $-factor $R$ containing $L$ and 
such that no component of $R$ induces a triangle, and
\\[0.5ex]
$(a3)$ if $L$ does not induce a triangle in $G$, 
then $G$ has  a $\Lambda $-factor, 
containing $L$ and a component that induces a triangle.
\es

\bp Let $L = xzz_1$.
Let $E'$ be the set of edges in $G$ that
belong to no triangle. 
Obviously, there is a natural bijection $\alpha : E(G') \to E'$. 
 Since each vertex of $G$ belongs to exactly one triangle,
we can assume that  $xz$ belongs to a triangle $T = xzs$.
\\[1ex]
${\bf (p1)}$ Suppose that $L$ induces a triangle in $G$, 
and so $s = z_1$.
Obviously the union of all triangles in $G$ contains
a $\Lambda $-factor, say $P$, of $G$ and 
$L \subset P$.
Therefore claim $(a1)$ is true.
\\[1ex]
${\bf (p2)}$ Now suppose that $L$ does not induce a triangle in $G$, and so $s \ne z_1$. 
Let $\bar{s} = ss_1$ and $\bar{z} = zz_1$ be the edges of $G$ not belonging to $T$, and therefore belonging to no triangles in $G$. Hence $\bar{s} = \alpha (\bar{s}')$ and 
$\bar{z} = \alpha (\bar{z}')$, where
$\bar{s}' = s's'_1$ and $\bar{z}' = z'z'_1$ are edges in $G'$, and $s'  = z'$.
Since every vertex in $G$ belongs to exactly one triangle, clearly $s_1 \ne z_1$.
\\[1ex]
${\bf (p2.1)}$
We prove $(a2)$. 
By using Tutte's criterion for a graph to have a perfect matching, it is easy to prove the following: 
\\[1ex]
{\sc Claim}.  
{\em If $A$ is a cubic 2-connected
graph, then for every 3-vertex path $J$ of $A$ there exists 
a 2-factor of $A$ containing $J$.}
 \\[1ex]
 \indent
By  the above {\sc Claim}, $G'$ has a 2-factor $F'$ containing 3-vertex path $S' = s'_1s'z'_1$.
Let $C'$ be the (cycle) component of $F'$ containing $S'$.
If $Q'$ is a (cycle) component of $F'$, then let $Q$ 
be the subgraph of $G$, induced by the edge subset 
$\{\alpha (e): e \in E(Q')\} \cup \{E(\Delta _v): v \in V(Q')\}$.
Obviously $v(Q) = 0 \bmod 3$ and $Q$ has a (unique) Hamiltonian cycle $H(Q)$. 
Also the union $F$ of all $Q$'s is a spanning subgraph of 
$G$ and each $Q$ is a component of $F$.  
Moreover, if $C$ is the component in $F$, corresponding to $C'$, then $L \subset H(C)$.
Therefore each $H(Q)$ has a  $\Lambda $-factor $P(Q)$, such that no component of $P(Q)$ induces a triangle, and 
$H(C)$ has a (unique) $\Lambda $-factor $P(C)$, such that $L \subset P(C)$ and no component of $P(C)$ induces 
a triangle.
The union of all these $\Lambda $-factors
is a $\Lambda $-factor $P$ of $G$ containing $L$ and 
such that no component of $P$ induces a triangle.
Therefore $(a2)$ holds.
\\[1ex]
${\bf (p2.2)}$
Now we prove $(a3)$.
Since $G'$ is 2-connected and cubic, there is a cycle $C'$ in $G'$ such that $V(C') \ne V(G')$ and $C'$ 
contains $S' = s'_1s'z'_1$.
Let, as above,  $C$ be the subgraph of $G$, induced by 
the edge subset
$\{\alpha (e): e \in E(C')\} \cup \{E(\Delta _v): v \in V(C')\}$.
Obviously, $v(C) = 0 \bmod 3$, $C$ has a (unique) Hamiltonian cycle $H$, and $L \subset H$.
Therefore $H$ has a (unique) $\Lambda $-factor $P(C)$ containing $L$.
Since $V(C') \ne V(G')$, we have $V(G' - C') \ne \emptyset $. Therefore $G - C$ has a triangle. Moreover, every vertex $v$ in $G - C$ belongs to a unique triangle $\Delta _v$, and therefore as in ${\bf (p1)}$, $G - C$ has a $\Lambda $-factor $Q$ whose every component induces a triangle in $G - C$. Then $P(C) \cup Q$ is a required a $\Lambda $-factor  in $G$.
\ep 
\\[2ex]
\indent
Theorem {\bf \ref{2con,cubic,tr}} is not true for a cubic, 2-connected, and claw-free graph $F$  with an edge $xy$ belonging to two triangles $T_i$ with $V(T_i) = \{x,y, z_i\}$ because $L = z_1xz_2$ is a 3-vertex path in $F$ and 
$y$ is an isolated vertex in $F - L$.

\bs 
\label{clfree,2con,(x,y)-}
Suppose that $G$ is a 2-connected  claw-free graph,
$v(G) = 2 \bmod 3$, and $x \in V(G)$. Then  there exist
at least two edges $xz_1$ and $xz_2$ in $G$  such that each $G - xz_i$ is connected and has a  $\Lambda $ -factor.
\es

\bp (uses {\bf \ref{2conclfr}}).
We need the following simple facts.
\\[1ex]
{\sc Claim} 1.
{\em Let $G$ be a 2-connected graph and,  
$x \in V(G)$. Then there exist
at least two edges $xs_1$ and $xs_2$ in $G$ such that each $G - xs_i$ is connected.}
\\[1ex]
{\sc Claim} 2.
{\em Let $G$ be a claw-free graph, $B$ is a 2-connected block of $G$, and $x$ is a boundary vertex  of $B$.
Then $B - x$ is either 2-connected or has exactly one edge.}
 \\[1ex]
\indent
By  {\sc Claim} 1, $G$ has an edge $xy$ such that 
$G - \{x,y\}$ is connected.
If $G - \{x, s\}$ for every $xs \in E(G)$,  then by {\sc Claim} 1, 
we are done.
Therefore we assume that $G - \{x, y\}$ is connected but has no $\Lambda$ - factor.

Then by {\bf \ref{2conclfr}}, $G -\{x,y\}$ has at least three end-blocks, say $B_i$, $i \in \{1, \ldots , k\}$, $k \ge 3$.
Let $b'_i$ be the boundary vertex of $B_i$. 
Let $V_i$ be the set of vertices in $\{x,y\}$ adjacent to the interior of $B_i$ and ${\cal B}_v$ be the set of the end-bocks in 
$G - \{x,y\}$ whose interior is adjacent to $v \in \{x,y\}$.
Since $G$ is 2-connected, each $|V_i| \ge 1$.
 Since $G$ is claw-free, each $|{\cal B}_v| \le 2$.
 Since $k \ge 3$, $|{\cal B}_v| = 2$ for some $v \in \{x,y\}$, say for $v = x$ and ${\cal B}_x = \{B_1,B_2\}$.
 Let $xb_i  \in E(G)$, where $b_i$ is an interior vertex of $B_i$, 
$i \in \{1,2\}$, and let $xb_j  \in E(G)$, where $b_j$ is an interior vertex of $B_j$, $j \ge 3$.
Since $G$ is claw-free, $\{x, y, b_1, b_2\}$
does not induce a claw in $G$. Therefore $yb_2 \in E(G)$.
If $k \ge 4$, then $\{y, b_2, b_3, b_4\}$
induces a claw in $G$, a contradiction.
Thus $k = 3$ and ${\cal B}_y = \{B_2, B_3\}$.

Suppose that $v(B_s) = 0 \bmod 3$ for some 
$s \in \{1,2,3\}$. Then $B_s$ is 2-connected.
Since $B_s$ is claw-free, by {\bf \ref{2conclfr}}, $B_s$ has a 
$\Lambda $-factor, say  $P$. 
Since $v(G) = 2 \bmod 3$, we have
$v(G - \{x,y, B_s\}) = 0 \bmod 3$.  
By {\sc Claim 2},  $G - \{x,y, B_s\}$ is claw-free, connected and has at most two end-blocks. 
Then by   {\bf \ref{2conclfr}},  $G - \{x,y, B_s\}$ has
$\Lambda $-factor, say $Q$. 
Therefore $P \cup Q$ is a $\Lambda $-factor of $G - \{x,y\}$, 
a contradiction.

Suppose that $v(B_r) = 1 \bmod 3$ for some 
$r \in \{1,2,3\}$. Then $B_r$ is 2-connected. 
Obviously $v(B_r -  b'_r) = 0 \bmod 3$ and claw-free. 
By {\sc Claim 2},  $B_r - b'_r$ is 2-connected. 
Then by {\bf \ref{2conclfr}}, $B_r$ has a 
$\Lambda $-factor, say  $P$. Obviously 
$G - \{x,y, B_r - b'_r\}$ is claw-free, connected and has at most two end-blocks. 
Then by   {\bf \ref{2conclfr}},  $G - \{x,y, B_s\}$ has
$\Lambda $-factor, say $Q$. 
Therefore $P \cup Q$ is a $\Lambda $-factor of $G - \{x,y\}$, 
a contradiction.

Now suppose that  $v(B_i) = 2 \bmod 3$ for every 
$i \in \{1,2,3\}$. 
By {\sc Claim 2}, either $B_i - \{b_i,b'_i\} $ is  
2-connected or 
$v(B_i - \{b_i,b'_i\}) = 0$ for every $i \in \{1,2,3\}$. 
In both cases by the arguments similar to that above, 
$G - \{x,b_i\}$ has a $\Lambda $-factor for $i \in \{1,2\}$ and
$G - \{y,b_i\}$ has a $\Lambda $-factor for $i \in \{2, 3\}$.
\ep
\\[2ex]
\indent
From {\bf \ref{clfree,2con,(x,y)-}} we have  for
3-connected claw-free graphs the following stronger result with a simpler proof.
\bs 
\label{clfree,3con,(x,y)-} 
Suppose that $G$ is a 3-connected  claw-free 
graph,  $v(G) = 2 \bmod 3$, and $xy \in E(G)$.
Then $G - \{x,y\}$ has a $\Lambda $-factor.
\es

\bp  (uses {\bf \ref{2conclfr}}). Let $G' = G - \{x,y\}$. Since $G$ is 3-connected, $G'$ is connected.
By {\bf \ref{2conclfr}}, it suffices to
prove that $G'$ has at most two end-blocks.
Suppose, on the contrary, that $G'$ has at least three end-blocks. 
Let $B_i$, $i \in \{1,2,3\}$, be some three blocks of $G'$.
Since $G$ is 3-connected, for every block $B_i$ and every vertex $v \in \{x,y\}$ there is an edge $vb_i$, where $b_i$ is an inner vertex of $B_i$.  Then $\{v, b_1, b_2, b_3\}$  induces a claw in $G$, a contradiction.
\ep
\\[2ex]
\indent
As we have seen in the proof of  {\bf \ref{clfree,2con,(x,y)-}}, the claim of {\bf \ref{clfree,3con,(x,y)-}} is not true for claw-free graphs of connectivity two.

\bs 
\label{clfree,3con,L-}
Suppose that $G$ is a 3-connected  claw-free graph,
$v(G) = 0 \bmod 3$,
and $xy \in E(G)$. Then  there exist
at least two 3-vertex paths $L_1$ and $L_2$ in $G$ centered at $y$, containing $xy$, and such that each $G - L_i$ is connected and has a  $\Lambda $-factor.
\es

\bp (uses {\bf \ref{2conclfr}}).
We need the following simple fact.
\\[1ex]
{\sc Claim} 1.
{\em Let $G$ be a 3-connected graph,  
$x \in V(G)$, and  $xy \in E(G)$. Then there exist
two 3-vertex paths $L_1$ and $L_2$ in $G$ centered at $y$, containing $xy$, and such that each $G - L_i$ is connected.}
\\[1ex]
\indent
By  {\sc Claim} 1, $G$ has a 3-vertex path $L = xyz$ such that
$G - L$ is connected. 
If every such 3-vertex path belongs to a $\Lambda $-factor of $G$, then by {\sc Claim} 1, 
we are done.
Therefore we assume that $G - L$ is connected but has no $\Lambda$ - factor.
Then by {\bf \ref{2conclfr}}, $G - L$ has at least three end-blocks, say $B_i$, $i \in \{1, \ldots , k\}$, $k \ge 3$.
Let $b'_i$ be the boundary vertex of $B_i$.
Let $V_i$ be the set of vertices in $L$ adjacent to inner 
vertices 
of $B_i$ and ${\cal B}_v$ be the set of the end-bocks in 
$G - L$ having an inner vertex 
adjacent to $v$ in $V(L)$.
Since $G$ is 3-connected, each $|V_i| \ge 2$.
 Since $G$ is claw-free, each $|{\cal B}_v| \le 2$.
It follows that $k = 3$, each $|V_i| = 2$, each $|{\cal B}_v| = 2$, as well as  all $V_i$'s are different and all ${\cal B}_v$'s are different. 
Let $s^1 = z$, $s^2 = x$,  $s^3 = y$, and 
$S = \{s^1, s^2, s^3\}$. 
We can assume that $V_i = S - s^i$, $i \in \{1,2,3\}$.
Then for every vertex $s^j \in V_i$ there is has a vertex  $b_i^j$ in  $B_i - b'_i$ adjacent to $s^j$, where $\{b_i^j: s^j \in V_i\}$ has exactly one vertex if and only if $B_i - b'_i$ has exactly one vertex.
Let $L_i = s^2s^3b_i$, where $b_i = b_i^3$.

By {\bf \ref{2conclfr}}, it surfices 
 to show that each $G - L_i$ is connected and has  at most two end-blocks. 

Let $i = 1$. If $B_1 - b_1$ is 2-connected, then 
$B_1 - b_1$ and $G - L_1 - (B_1 - b'_1)$ are the two 
end-blocks of $G - L_1$ and we are done. If $B_1 - b_1$ is empty, then $G - L_1$ is 2-connected.
So we assume that
$B_1 - b_1$ is not empty and not 2-connected. Then $B_1 - b_1$ is connected and has exactly two end-blocks, say $C_1$ and $C_2$. Let $c'_i$ be the boundary vertex of $C_i$ in $B_1 - b_1$. Since $G$ is 3-connected, each $C_i -c'_i$ has a vertex  adjacent to $\{s^2,s^3\}$.
We can assume that a vertex $c_1$ in $C_1 - c'_1$ is adjacent to $s^2$.
If there exists a vertex $c_2$  in $C_2 - c'_2$ adjacent to 
$s^2$,
then 
$\{s^2, b_3^2, c_1, c_2\}$ 
induces a claw in $G$, a contradiction. So suppose that  no vertex in $C_2 - c'_2$ is adjacent to $s^2$.
Then there is a vertex $c_2$ in $C_2 - c'_2$ adjacent to 
$s^3$.
Then $\{s^2, s^3, b_2^3, c_2\}$
 induces a claw in $G$, a contradiction.

Now let $i = 2$.
If $B_2 - b_2$ is 2-connected, then 
$B_1$ and $G - L_2 - (B_1 - b'_1)$ are the two 
end-blocks of $G - L_2$ and we are done. If $B_1 - b_1$ is empty, then $G - L_2$ has two end-blocks, namely $B_1$ and the subgraph of $G$ induced by $B_3 \cup s^1$.
So we assume that
$B_2 - b_2$ is not empty and not 2-connected. Then $B_2 - b_2$ is connected and has exactly two end-blocks, say $D_1$ and $D_2$. Let $d'_i$ be the boundary vertex of $D_i$ in $B_2 - b_2$. Since $G$ is 3-connected, each $D_i -d'_i$ has a  vertex  adjacent to $\{s^1,s^3\}$.
We can assume that a vertex $d_1$ in $D_1 - d'_1$ is adjacent to $s^3$.
If there exists a vertex $d_2$  in $D_2 - d'_2$ adjacent to 
$s^3$,  
then $\{s^3,d_1,d_2, b_1^3\}$ induces a claw in $G$, a contradiction.
So suppose that no vertex in $D_2 - d'_2$ is adjacent to 
$s^3$.
Then there is a vertex $d_2$ in $D_2 - d'_2$ adjacent to 
$s^1$.
Then $\{s^1, s^3, b_3^1, d_2\}$  
induces a claw in $G$, a contradiction.
\ep
\\[2ex]
\indent
From the proof of {\bf \ref{clfree,3con,L-}} we  have, in particular:
\bs 
\label{clfree,3con,deg3cntrL-}
Suppose that $G$ is a 3-connected  claw-free graph.
If $L$ is a 3-vertex path and the center vertex of $L$ has degree 3 in $G$, then  $G - L$ is connected and has a  $\Lambda $-factor in $G$.
\es

Obviously, {\bf \ref{2con,cubic,tr}} (a) follows from 
{\bf \ref{clfree,3con,deg3cntrL-}}.
\\[2ex]
\indent
From the proof of {\bf \ref{clfree,3con,L-}} we also have:
\bs 
\label{clfree,4con,deg3cntrL-}
Suppose that $G$ is a 4-connected  claw-free graph.
Then $G - L$ is connected and has a  $\Lambda $-factor for every 3-vertex path $L$ in $G$.
\es

The claim of {\bf \ref{clfree,4con,deg3cntrL-}} may not be true for a  claw-free graph of connectivity 3 if they are not cubic.
A graph obtained  obtained from a claw  by replacing its vertex of degree 3 by a triangle is called a {\em net}. 
Let $N$ be a net with the three leaves $v_1$, $v_2$, and $v_3$,
$T$  a triangle with $V(T) = \{t_1, t_2, t_3\}$, and let $N$ and $T$ be disjoint. 
Let $H = N \cup T \cup \{v_it_j: i, j \in \{1,2,3\}, i \ne j\}$.
Then $H$ is a 3-connected claw-free graph, $v(H) = 9$, each $d(t_i, H) = 4$, $d(x, H) = 3$ for every $x \in V(H - T)$, and  
$ H - T = N$ has no $\Lambda $-factor.  If $L$ is a 3-vertex path in $T$, then 
$H - L = H - T$, and so $H - L$ has no $\Lambda $-factor.
There are infinitely many pairs $(G, L)$ such that
$G$ is a 3-connected, claw-free, and non-cubic graph, 
$v(G) = 0 \bmod 3$, $L$ is a 3-vertex path in $G$, and 
$G - L$ has no 
$\Lambda $-factor. By {\bf \ref{A}}, such a pair can be obtained from the above pair $(H,L)$ by replacing $N$ by any graph $A$ with three leaves from the class ${\cal A}$ (defined below before {\bf \ref{A}}) provided $v(A) = 0 \bmod 3$.
\\[2ex]
\indent
From {\bf \ref{clfree,3con,L-}} we have, in particular:
\bs 
\label{clfree,3con,einL}
Suppose that $G$ is a 3-connected  claw-free graph and 
$e \in E(G)$.
Then 
\\[.5ex]
$(a1)$ there exists a $\Lambda $-factor in $G$ containing $e$ and 
\\[.5ex]
$(a2)$ there exists a $\Lambda $-factor in $G$ avoiding $e$, i.e. $G - e$ has a $\Lambda $-factor.
\es

The following examples show that condition 
``{\em $G$ is a 3-connected graph}'' in 
{\bf \ref{clfree,3con,einL}} is essential for claim $(a1)$.
Let $R$ be the graph obtained from two disjoint cycles $A$ and $B$ by adding a new vertex $z$,  and the set of new edges $\{a_iz, b_iz: i \in \{1,2\}\}$, 
where $a = a_1a_2 \in E(A)$ and $b = b_1b_2 \in E(B)$. 
It is easy to see that $Q$ is a claw-free graph of connectivity one. 
Furthermore, if $v(A) = 1 \bmod 3$ and 
$v(B) = 1 \bmod 3$, then $v(Q) = 0 \bmod 3$ and 
$Q$ has no 
$\Lambda $-factor containing edge $e \in \{a, b\}$.
Similarly,
let $Q$ be the graph obtained from two disjoint cycles $A$ and $B$ by adding two new vertices $z_1$ and $z_2$, a new edge $e = z_1z_2$,  and the set of new edges $\{a_iz_j, b_iz_j: i,j \in \{1,2\}\}$, where $a_1a_2 \in E(A)$ and $b_1b_2 \in E(B)$. 
It is easy to see that $Q$ is a claw-free graph of connectivity two. Furthermore, if $v(A) = 2 \bmod 3$ and 
$v(B) = 2 \bmod 3$, then $v(Q) = 0 \bmod 3$ and $Q$ has no 
$\Lambda $-factor containing edge $e$.

As to claim $(a2)$ in {\bf \ref{clfree,3con,einL}}, it turns out that this claim is also true for 2-connected claw-free graphs.
\bs 
\label{clfree,2con-avoid-e}
Suppose that $G$ is a 2-connected  claw-free graph, 
$v(G) = 0 \bmod 3$, and $e \in E(G)$.
Then  $G - e$ has a $\Lambda $-factor.
\es

\bp 
A graph $H$ is called {\em minimal 2-connected} if $H$ is 2-connected  but $H - u$ is not 2-connected for every $u \in E(H)$. A {\em frame} of a graph $G$ is a minimal 2-connected spanning subgraph of $G$. Clearly, every 
2-connected graph has a frame. In \cite{KKN} we describe  Procedure 1 that provides an ear-assembly $A$ of a special frame  of 
a 2-connected claw-free graph. In particular, the last ear  of $A$ contains a $\Lambda $-packing $P$ such that 
$G - P$ is also 2-connected claw-free graph. We modify Procedure 1 by replacing the first step of this procedure 
``{\em Find a longest cycle $G_0$ in $G$}'' by ``{\em Find a longest cycle  $G'_0$ among all cycles $C$ in $G$ such that  edge $e$  either belongs to $C$ or is a chord of $C$}''. 
Since $G$ is 2-connected, $G$ has a cycle containing $e$. Therefore a cycle $G'_0$ does exist.
Then the resulting Procedure ${\cal P}$ provides  
an ear-assembly of a frame  of $G$ with the property  that the last ear of this frame has a $\Lambda $-packing $Q$ such that $e \not \in E( Q)$ and $G - Q$ is a 2-connected claw-free graph that may contain $e$.

We can use  Procedure ${\cal P}$ to prove our claim by induction on $v(G)$. If $G$ is a cycle containing $e$, then our claim is  obviously true. Procedure  ${\cal P}$ mentioned above guarantees the existence of a $\Lambda $-packing $Q$ such that $Q$ avoids $e$
and $G - Q$ is a 2-connected claw-free graph that may contain $e$. Obviously $v(G - Q) = 0 \bmod 3$ and 
$v(G - Q) < v(G)$. By the induction hypothesis, $G - Q$ has a $\Lambda $-factor $R$ avoiding $e$. Then 
$Q \cup R$ is a $\Lambda $-factor of $G$ avoiding edge  $e$.
\ep
\\[2ex]
\indent
We need the following fact interesting in itself.
Let ${\cal A}$ denote the set of graphs $A$ with the following properties:
\\[.5ex]
$(c1)$ $A$ is connected,
\\[.5ex]
$(c2)$ every vertex in $A$ has degree at most 3,
\\[.5ex]
$(c3)$ every vertex in $A$ of degree 2 or 3 belongs to exactly one triangle, 
and
\\[.5ex]
$(c4)$ $A$ has exactly three vertices of degree 1 which we call the {\em leaves} of $A$.
\bs
\label{A}
 If $A \in {\cal A}$, then $A$ has no $\Lambda $-factor.
 \es

\bp Let $A \in {\cal A}$. If $v(A) \ne 0 \bmod 3$, then our claim is clearly true. So we assume that  $v(A) = 0\bmod 3$.
We prove our claim by induction on $v(G)$. The smallest graph in ${\cal A}$ is a net $N$ with $v(N) = 6$ and our claim is obviously true for $N$. So let $v(A) \ge 9$.
Suppose, on the contrary, that $A$ has a $\Lambda $-factor $P$. Let $v$ be a leaf of
$A$ and $vx$ the edge incident to $v$. 
Since $P$ is a  $\Lambda $-factor in $A$, it has a component $L = vxy$, and so $P - L$ is a $\Lambda $-factor in $A - L$ and $d(x, A) \ge 2$.
By property $(c3)$, $x$ belongs to a unique triangle $xyz$ in $A$ and $d(x,a) = 3$, and so $s \in \{y,z\}$.  
If $d(z,A) = 2$, then $z$ is an isolated vertex in $A - L$, and so $P$ is not a  $\Lambda $-factor in $A$, a contradiction.
Therefore by $(c2)$, $d(z,A) = 3$.
Therefore $A - L$ satisfies $(c2)$, $(c3)$, and $(c4)$.

Suppose that $G - L$ is not connected and that the three leaves 
do not belong to a common component. Then $A - L$ has a component $C$ with $v(C) \ne 0 \bmod 3$, and so $A - L$ has no $\Lambda $-factor, a contradiction.

Now suppose that $A - L$ has a component $C$ containing all three leaves of $A - L$. Then $C \in {\cal A}$ and 
$v(C) < v(A)$. By the induction hypothesis, $C$ has no $\Lambda $-factor. Therefore $A - L$ also has no $\Lambda $-factor, a contradiction.
\ep
\\[2ex]
\indent
Given $E \subseteq E(G)$, let $\dot{E}$ denote the subgraph of $G$ induced by $E$.

\bs 
\label{clfree,2con-avoid(a,b,c)}
Suppose that $G$ is a cubic 2-connected graph and  that every vertex in $G$ belongs to exactly one  triangle  
$($and so $G$ is claw-free$)$, i.e. $G = F^\Delta $, where $F$ is a cubic 2-connected graph.
Let  $E \subset E(G)$ and  $|E| = 3$. Then 
 the following are equivalent:
\\[.5ex]
$(g)$
$G - E$ has no $\Lambda $-factor and
\\[.5ex]
$(e)$ $\dot{E}$ satisfies one of the following conditions:
\\[.5ex]
\indent
$(e1)$ $\dot{E}$ is  a claw,
\\[.5ex]
\indent
$(e2)$ $\dot{E}$ is a triangle,
\\[.5ex]
\indent
$(e3)$
$\dot{E}$ has exactly two components, 
the 2-edge component $\dot{E}_2$  belongs to a triangle in $G$,
the 1-edge component $\dot{E}_1$ belongs to no triangle in $G$, and $G - E$ is not connected, and 
\\[.5ex]
\indent
$(e4)$
$\dot{E}$ has exactly two components,  
in $G$ the 2-edge component  $\dot{E}_2$ belongs to a triangle, say $T$, the 1-edge component $\dot{E}_1$ also belongs to a triangle, say $D$, and $\dot{E}_1$, $\dot{E}_2$ belong to different component of $G - \{d,t\}$, where
$d$ is the edge incident to the vertex of $D - \dot{E}_1$ and 
$t$ is the edge  in $G - E$ incident to the isolated vertex of 
$T - E$.
 \es

\bp (uses 
{\bf \ref{2conclfr}}, 
{\bf \ref{2con,cubic,tr}}(a), and
{\bf \ref{A}}). 
Let $X, Y \subset E(G)$ such that $X$ meets no triangle in $G$, each edge in $Y$ belongs to a triangle in $G$,  and  no triangle in $G$ has more than one edge from $Y$, and so 
$X \cap Y = \emptyset $.
We will use the following simple observation.
\\[.5ex]
{\sc Claim.}
{\em 
$G - X - Y$ has a  $\Lambda $-factor $P$ such that  every component of $P$ induces a triangle in $G$ and if an edge $y$ from $Y$ is in a triangle $T$, then $T - y$ is a component of $P$.}
\\[.5ex]
\indent
Let $E = \{a,b,c\}$. By the above {\sc Claim}, we can assume that edges $a$ and $b$ belong to the same triangle $T$. \\[1ex]
${\bf (p1)}$ We prove $(e) \Rightarrow (g)$.

Suppose that $\dot{E}$ satisfies $(e1)$, i.e. $\dot{E}$ is a claw.
Then $G - E$ has an isolated vertex and therefore has  no $\Lambda $-factor. 

Suppose that $\dot{E}$ satisfies $(e2)$, i.e. $\dot{E}$ is a triangle. 
Then $G - E \in {\cal A}$. By {\bf \ref{A}}, $G - E$ has no 
$\Lambda $-factor.

Suppose that $\dot{E}$ satisfies $(e3)$, i.e. $G - E$ is not connected and
$\dot{E}$ has exactly two components $\dot{E}_2$ and $\dot{E}_1$
induced by $\{a,b\}$ and $c$, respectively, where 
$\dot{E}_2$ belongs to the triangle $T$ but 
$\dot{E}_1$ belongs to no triangle in $G$.
Then $t$ is the dangling edge in $G - E$. 
Let  
$S$ be the component in $G - E$ containing edge $t$. 
Then every vertex in $S$ distinct from the leaf incident to $t$ belongs to exactly one triangle. Therefore $v(S) = 1\bmod 3$. 
Thus $G - E$ has no $\Lambda $-factor. 

Now suppose that $\dot{E}$ satisfies $(e4)$.
By $(e4)$, $t$ is the edge in $G - E$ incident to $z$. 
Suppose, on the contrary, that $G - E$ has a 
$\Lambda $-factor, say $P$.
Since $G - \{d,t\}$ is not connected, $G - E - E(D)$ is also not connected.
Obviously the component of $G - E - E(D)$ containing $z$ belongs to ${\cal A}$. 
Therefore by {\bf \ref{A}}, $G - E - E(D)$ has no  $\Lambda $-factor. Thus $P$ has a 3-vertex path $L$ containing exactly one edge in $D$ adjacent  to edge $d$.
Now if $C$ is a component of $G - E - L$, then 
$v(C) \ne 0 \bmod 3$. 
Therefore  $G - E$ has no $\Lambda $-factor, a contradiction.
\\[1ex]
${\bf (p2)}$ Now we prove $(g) \Rightarrow (e)$.
Namely, we assume that $\dot {E}$ does not satisfy $(e)$ and we want to show that in this case $G - E$ has a 
$\Lambda $-factor. Let $u$ be the edge of $T$ distinct from $a$ and $b$.

Suppose that $\dot {E}$ is connected, and so $\dot {E}$ is a 3-edge path.   
Let $V$ be a 3-vertex path in $G$ containing $u$ and avoiding $E$. 
Then $G - V$ has no edges from $E$,  and so 
$G - V = G - E -V$.
By {\bf \ref{2con,cubic,tr}}(a), $G - V$ has a $\Lambda $-factor.

Now suppose that $\dot {E}$ is not connected, and so 
$\dot {E}$ has exactly two components induced by 
$\{a,b\}$ and by $c$, respectively.
Since $\dot {E}$ does not satisfy $(e)$, $\dot {E}$ is not a claw and not a triangle, and so $u, t \not \in E$.
\\[1ex]
${\bf (p2.1)}$ Suppose that $c$ belongs to no triangle in $G$.
Since $\dot {E}$ does not satisfy $(e)$, $G - E$ is connected. Clearly, $G - E$ is claw-free. 
Also $G - E$ has exactly two end-blocks and the  block of one edge $t$ is one of them.
By {\bf \ref{2conclfr}}, $G - E$ has a $\Lambda $-factor.
\\[1ex]
${\bf (p2.2)}$ Now suppose that $c$ belongs to a triangle $D$ in $G$. Then $D \ne T$.
Let $G' = G - D - \{a,b\}$. Then $G'$ is claw-free and has no edges from $E$.

Suppose that $G'$ is connected. Then as in ${\bf (p2.1)}$,
$G'$ has exactly two end-blocks. By {\bf \ref{2conclfr}}, $G'$ has a $\Lambda $-factor, say $P$. 
Let $L = D - c$, and so $L$ is a 3-vertex path.
Then $P \cup \{L\}$ is a $\Lambda $-factor in $G - E$.

Now suppose that $G'$ is not connected. Let $C$ be the component of $G'$ containing edge $t$ and $q$ the edge connecting $C$ and $D$.
Let $L$ be a 3-vertex path in $G$ containing $q$ and an edge in $D - c$. It is sufficient to show that 
$G - E - L$ has a $\Lambda $-factor. Obviously  $G - E - L$ is claw-free. Let $Q$ be a component of $G - E - L$. 
Since $\dot {E}$ does not satisfy $(e4)$, 
$v(Q) = 0\bmod 3$ and $Q$ has exactly two end-blocks. 
By  {\bf \ref{2conclfr}}, $Q$ has a $\Lambda $-factor.
Therefore $G - E - L$ also has a $\Lambda $-factor.
\ep
\\[2ex]
\indent
From {\bf \ref{clfree,2con-avoid(a,b,c)}} we have, in particular:
\bs 
\label{clfree,2con-avoid(a,b)}
Suppose that  $G = F^\Delta $, where $F$ is a cubic
 2-connected graph $($and so $G$ is claw-free$)$.
Let  $E \subset E(G)$ and  $|E| = 2$. Then
$G - E$ has a $\Lambda $-factor.
\es

From {\bf \ref{clfree,2con-avoid(a,b,c)}} we also have:

\bs 
\label{clfree,3con-avoid(a,b,c)}
Suppose that  $G$ is a 3-connected claw-free graph.
Let  $E \subset E(G)$ and  $|E| = 3$. Then
$G - E$ has a $\Lambda $-factor if and only if 
$\dot {E}$ is not a claw and not a traingle.
\es

\bs 
\label{clfree,2con,x-}
Suppose that $G$ is a 2-connected  claw-free graph, 
$v(G) = 1 \bmod 3$  and  $x \in V(G)$.
Then $G - x$ has a $\Lambda $-factor. 
\es

\bp (uses {\bf \ref{2conclfr}}).
Let $x \in V(G)$. 
Since $v(G) = 1 \bmod 3$, clearly $v(G - x) = 0 \bmod 3$.
Since $G$ is 2-connected, $G - x$ is connected.
Since $G$ is  claw-free, $G - x$ is claw-free and has 
at most two end-blocks.
By {\bf \ref{2conclfr}},
$G - x$ has a $\Lambda $-factor.
\ep

\bs 
\label{clfree,3con,(x,e)-} 
Suppose that $G$ is a 3-connected  claw-free 
graph,  $v(G) = 1 \bmod 3$, $x \in V(G)$ and $e \in E(G)$.
Then $G - \{x,e\}$ has a $\Lambda $-factor.
\es

\bp (uses {\bf \ref{clfree,2con-avoid-e}} and 
{\bf \ref{clfree,2con,x-}}).
Since $G$ is 3-connected, $G - x$ is a 2-connected  claw-free graph. Since $v(G) = 1 \bmod 3$, we have
$v(G - x) = 0 \bmod 3$. 
By {\bf \ref{clfree,2con,x-}}, $G - x$ has a $\Lambda $-factor $P$.
If $e \not \in E(G - x)$, then $P$ is a $\Lambda $-factor  of 
$G - \{x,e\}$. If $e  \in E(G - x)$, then by 
{\bf \ref{clfree,2con-avoid-e}}, $G - \{x,e\}$ has 
a $\Lambda $-factor.
\ep
\\[2ex]
\indent
Obviously, the claim in {\bf \ref{clfree,3con,(x,e)-}} may not be true  for a  claw-free graph of connectivity 2.
\\[2ex]
\indent
From  {\bf \ref{2con,cubic,tr}},  {\bf \ref{clfree,3con,(x,y)-}}, 
{\bf \ref{clfree,3con,einL}}, {\bf \ref{clfree,2con-avoid(a,b)}},
{\bf \ref{clfree,2con,x-}}, and {\bf \ref{clfree,3con,(x,e)-}} 
we have, in particular:
\bs 
\label{clfree,3con,cubic,ztf}
All claims in  {\bf \ref{cubic3-con}} except for $(z5)$ are true for 3-connected claw-free graphs and $(z5)$ is true for 
cubic, 2-connected graphs such that every vertex belongs to exactly one triangle. 
\es


\end{document}